\documentclass [12pt]{article}
\usepackage{graphicx,amssymb,amsfonts,latexsym,amsmath,amsthm,times}
\usepackage{epsfig}
\usepackage{color}
\usepackage[a4paper]{geometry}
\def\lanbox{\hbox{$\, \vrule height 0.25cm width 0.25cm depth 0.01cm \,$}}

\numberwithin{equation}{section}

\begin{document}

\vspace*{1.4cm}

\normalsize \centerline{\Large \bf EXISTENCE OF STATIONARY SOLUTIONS
FOR}

\medskip

\centerline{\Large\bf SOME INTEGRO-DIFFERENTIAL EQUATIONS WITH}

\medskip

\centerline{\Large\bf THE DOUBLE SCALE ANOMALOUS DIFFUSION}

\vspace*{1cm}

\centerline{\bf Vitali Vougalter$^{1 \ *}$, Vitaly Volpert$^{2,3}$}

\vspace*{0.5cm}

\centerline{$^{1 \ *}$ Department of Mathematics, University
of Toronto}

\centerline{Toronto, Ontario, M5S 2E4, Canada}

\centerline{ e-mail: vitali@math.toronto.edu}

\medskip

\centerline{$^2$ Institute Camille Jordan, UMR 5208 CNRS,
University Lyon 1}

\centerline{ Villeurbanne, 69622, France}

\centerline{$^3$ Peoples' Friendship University of Russia, 6 Miklukho-Maklaya
St,}

\centerline{Moscow, 117198, Russia}

\centerline{e-mail: volpert@math.univ-lyon1.fr}

\medskip


\vspace*{0.25cm}

\noindent {\bf Abstract:}
The paper is devoted to the investigation of the solvability of an
integro-differential equation in the case of the double scale anomalous
diffusion with a sum of two negative Laplacians in different fractional powers in ${\mathbb R}^{3}$. The proof
of the existence of solutions relies on a f{i}xed point technique. Solvability
conditions for the elliptic operators without the Fredholm property in
unbounded domains are used.

\vspace*{0.25cm}

\noindent {\bf AMS Subject Classi{f}ication:} 35R11, 35P30, 45K05


\vspace*{0.5cm}

\bigskip

\bigskip


\setcounter{equation}{0}

\section{\bf Introduction}

The present article is devoted to the studies of the existence of the stationary
solutions of the following nonlocal integro-differential problem
\begin{equation}
\label{h}
\frac{\partial u}{\partial t} =
-D[(-\Delta)^{s_{1}}+ (-\Delta)^{s_{2}}]u +
\int_{{\mathbb R}^{3}}K(x-y)g(u(y,t))dy + f(x),
\end{equation}
where $\displaystyle{\frac{1}{4}<s_{1}<\frac{3}{4}}$ and
$\displaystyle{s_{1}<s_{2}<1}$. The equations of this kind
are crucial to the cell population dynamics. The results of the work
are obtained in these particular ranges of the values of the parameters
$s_{1}$ and $s_{2}$ in
the powers of the negative Laplacians. This is based on the
solvability of the linear Poisson type equation (\ref{lp}) and the
applicability of the Sobolev inequality for the fractional Laplace operator
(\ref{frs}). The solvability of the problem analogous to (\ref{h}) containing
a single fractional Laplacian in the diffusion term was considered in
~\cite{VV15}.
The space variable $x$ here
corresponds to the cell genotype,
$u(x,t)$ stands for the cell density as a function of their genotype and time.
The right side of our equation describes the evolution of the cell density by
means of the cell proliferation, mutations and cell influx/efflux.
The double scale anomalous diffusion term in this context describes
the change of genotype due to the small random mutations, and the integral
production term corresponds to large
mutations. The function $g(u)$ designates the rate of the cell birth depending
on $u$
(density dependent proliferation), and the kernel $K(x-y)$ stands for
the proportion of the newly born cells changing their genotype from
$y$ to $x$.
We assumed here that it depends on the distance between the genotypes.
The last term in the right side of (\ref{h}) denotes
the influx or efflux of cells for different genotypes.

The fractional Laplace operator describes a particular case
of the anomalous diffusion actively used in the context of different
applications in plasma physics and
turbulence \cite{2}, \cite{1}, surface diffusion \cite{3}, \cite{4},
semiconductors \cite{5} and so on. The anomalous
diffusion can be described as a random process of the particle motion
characterized by the probability density distribution of the jump length.
The moments of this density distribution are finite in
the case of the normal diffusion, but this is not the case for the anomalous
diffusion. The asymptotic behavior at the infinity of the probability density
function
determines the value of the power of the negative Laplace operator
(see ~\cite{MK00}).
In the current work we will treat the situation when
$\displaystyle{\frac{1}{4}<s_{1}<\frac{3}{4}, \ s_{1}<s_{2}<1}$ in the
space of three dimensions. 

The solvability of the integro-differential
equation in the case of the standard
Laplace operator in the diffusion term was considered in ~\cite{VV111}.
A Faber-Krahn inequality for mixed local and nonlocal operators was proved
in ~\cite{BDVV23}.
Existence and dynamics of normalized solutions to nonlinear Schr\"odinger
equations with mixed fractional Laplacians were investigated in ~\cite{CGH23}.
An existence theory for superposition operators of mixed order subject to
jumping nonlinearities was developed in ~\cite{DPSV24}.
The necessary conditions of the preservation of the nonnegativity of the
solutions of a system of parabolic equations in the situation of the double
scale anomalous diffusion were obtained in ~\cite{EV221}. The work
~\cite{GNTUV21} is devoted to the simultaneous inversion for the fractional
exponents in the space-time fractional diffusion problem. A discretization
method for nonlocal diffusion type equations was developed in ~\cite{MORV22}.

We set $D=1$ and demonstrate the existence of solutions of the equation
\begin{equation}
\label{p}
-[(-\Delta)^{s_{1}}+(-\Delta)^{s_{2}}]u +\int_{{\mathbb R}^{3}} K(x-y)g(u(y))dy +
f(x) = 0
\end{equation}
with
$\displaystyle{\frac{1}{4}<s_{1}<\frac{3}{4}, \ s_{1}<s_{2}<1}$.
Let us discuss the case when the linear part of such operator does not
satisfy the Fredholm property. Consequently,  the conventional methods of
the nonlinear analysis may not be applicable. We use the solvability conditions
for the non-Fredholm operators along with the method of contraction mappings.

Consider the problem
\begin{equation}
\label{eq1}
 -\Delta u + V(x) u - a u=f,
\end{equation}
where $u \in E= H^{2}({\mathbb R}^{d})$ and  $f \in F=
L^{2}({\mathbb R}^{d}), \ d\in {\mathbb N}$, $a$ is a constant and
the scalar potential function $V(x)$ is either trivial
or converges to $0$ at infinity. For $a \geq 0$, the essential spectrum of the
operator $A : E \to F$, which corresponds to the left  side of equation
(\ref{eq1}) contains the origin. As a consequence, this operator fails to
satisfy the Fredholm property. Its image is not closed, for $d>1$
the dimension of its kernel and the codimension of its image are
not finite. The present article deals with the studies of some properties
of the operators of this kind. Note that the elliptic equations with
non-Fredholm operators were considered actively in recent years.
Methods in weighted Sobolev and H\"older spaces were developed in
~\cite{Amrouche1997}, ~\cite{Amrouche2008}, ~\cite{Bolley1993},
~\cite{Bolley2001}, ~\cite{B88}. The non-Fredholm Schr\"odinger type operators
were approached with the methods of the spectral and the
scattering theory in ~\cite{V2011}, ~\cite{VV08}.
The nonlinear non-Fredholm elliptic problems were considered in ~\cite{EV22},
~\cite{EV221}, ~\cite{VV111}, ~\cite{VV14}, ~\cite{VV15}.
The interesting applications to the theory of the reaction-diffusion
problems were developed in ~\cite{DMV05}, ~\cite{DMV08}. Fredholm structures,
topological invariants and applications were treated in ~\cite{E09}.
The works ~\cite{GS05} and ~\cite{RS01} are important for the understanding
of the Fredholm and properness properties of the quasilinear elliptic systems
of the second order and of the operators of this kind on ${\mathbb R}^{N}$.
The operators without
the Fredholm property appear also when studying the wave systems with an
infinite
number of localized traveling waves (see ~\cite{AMP14}). The standing lattice
solitons in the discrete NLS equation with saturation were discussed in
~\cite{AKLP19}.
In particular,
when $a=0$ the operator $A$ is Fredholm in some properly chosen weighted
spaces (see \cite{Amrouche1997}, \cite{Amrouche2008}, \cite{Bolley1993},
\cite{Bolley2001}, \cite{B88}). However, the situation of $a \neq 0$ is
significantly
different and the method developed in these articles cannot be applied. The
front
propagation equations with the anomalous diffusion were considered actively in
recent years (see e.g. ~\cite{VNN10}, ~\cite{VNN13}).

\medskip

Let us set $K(x) = \varepsilon {\cal K}(x)$ with $\varepsilon \geq 0$
and assume that the conditions below are satis{f}ied.

\medskip

\noindent
{\bf Assumption 1.1.}  {\it Let $\displaystyle{\frac{1}{4}<s_{1}<\frac{3}{4}}$
and
$\displaystyle{s_{1}<s_{2}<1}$.
Suppose that $f(x): {\mathbb R}^{3}\to {\mathbb R}$ is nontrivial, such that
$f(x)\in L^{1}({\mathbb R}^{3})$ and
$(-\Delta)^{1-s_{1}} f(x)\in
L^{2}({\mathbb R}^{3})$. Assume also that
${\cal K}(x): {\mathbb R}^{3}\to {\mathbb R}$ and
${\cal K}(x)\in L^{1}({\mathbb R}^{3})$. In addition,
$(-\Delta)^{1-s_{1}}{\cal K}(x)\in L^{2}({\mathbb R}^{3})$,
so that}
$$
Q:=\big\|(-\Delta)^{1-s_{1}}{\cal K}(x)\big\|_{L^{2}({\mathbb R}^{3})}>0.
$$

\bigskip

We choose the space dimensions $d=3$. This is related to the solvability
conditions for the linear Poisson type problem (\ref{lp}) formulated in
Lemma 4.1 below. From the perspective of the applications, the space
dimension is not limited to $d=3$ because the space variable corresponds to
the cell genotype but not to the usual physical space.

Let us apply the Sobolev
inequality for the fractional Laplace operator
(see Lemma 2.2 of ~\cite{HYZ12}, also ~\cite{L83}), namely
\begin{equation}
\label{frs}
\|f(x)\|_{L^{\frac{6}{4s_{1}-1}}({\mathbb R}^{3})}\leq c_{s_{1}}
\|(-\Delta)^{1-s_{1}}f(x)\|_
{L^{2}({\mathbb R}^{3})}, \quad \frac{1}{4}<1-s_{1}<\frac{3}{4}
\end{equation}
along with Assumption 1.1 above and the standard interpolation argument.
This gives us
\begin{equation}
\label{fl2}
f(x)\in  L^{2}({\mathbb R}^{3})
\end{equation}
as well.
For the technical purposes, we introduce the Sobolev space
$$
H^{2s_{2}}({\mathbb R}^{3}):=\big\{u(x):{\mathbb R}^{3}\to {\mathbb R} \ | \
u(x)\in L^{2}({\mathbb R}^{3}), \ (-\Delta)^{s_{2}} u \in L^{2}
({\mathbb R}^{3}) \big\},
$$
where $\displaystyle{0<s_{2}\leq 1}$, which is
equipped with the norm
\begin{equation}
\label{n}
\|u\|_{H^{2s_{2}}({\mathbb R}^{3})}^{2}:=\|u\|_{L^{2}({\mathbb R}^{3})}^{2}+
\big\|(-\Delta)^{s_{2}}u\big\|_{L^{2}({\mathbb R}^{3})}^{2}.
\end{equation}
The standard Sobolev embedding yields that
\begin{equation}
\label{e}
\|u\|_{L^{\infty}({\mathbb R}^{3})}\leq c_{e}\|u\|_{H^{2}({\mathbb R}^{3})},
\end{equation}
where $c_{e}>0$ is the constant of the embedding. Here
\begin{equation}
\label{h3}
\|u\|_{H^{2}({\mathbb R}^{3})}^{2}:=\|u\|_{L^{2}({\mathbb R}^{3})}^{2}+
\big\|\Delta u\big\|_{L^{2}({\mathbb R}^{3})}^{2}.
\end{equation}
When the nonnegative parameter $\varepsilon$ vanishes, we arrive at the
linear Poisson type equation (\ref{lp}).
By means of Lemma 4.1 below along with Assumption 1.1, equation (\ref{lp})
admits a unique solution
$$
u_{0}(x)\in H^{2s_{2}}({\mathbb R}^{3}), \quad \frac{1}{4}<s_{1}<\frac{3}{4}, \quad
s_{1}<s_{2}<1,
$$
such that no orthogonality conditions are required.
According to Assumption 1.1,
\begin{equation}
\label{s1s2u0f}
[-\Delta+
(-\Delta)^{1+s_{2}-s_{1}}] u_{0}(x)=(-\Delta)^{1-s_{1}}f(x)\in
L^{2}({\mathbb R}^{3}).
\end{equation}
It follows easily from (\ref{s1s2u0f}) via the standard Fourier
transform (\ref{f}) that $\Delta u_{0}(x)\in L^{2}({\mathbb R}^{3})$.
By virtue of the definition of the norm (\ref{h3}),
we derive for the unique solution of linear problem (\ref{lp}) that
$u_{0}(x)\in H^{2}({\mathbb R}^{3})$. Note that $u_{0}(x)$ is nontrivial
since the function $f(x)$ does not vanish identically in ${\mathbb R}^{3}$
as assumed.

Let us seek the resulting solution of nonlinear equation (\ref{p}) as
\begin{equation}
\label{r}
u(x)=u_{0}(x)+u_{p}(x).
\end{equation}
Evidently, we obtain the perturbative problem
\begin{equation}
\label{pert}
[(-\Delta)^{s_{1}}+(-\Delta)^{s_{2}}]u_{p}(x)=\varepsilon \int_{{\mathbb R}^{3}}
{\cal K}(x-y)g(u_{0}(y)+u_{p}(y))dy,
\end{equation}
where
$\displaystyle{\frac{1}{4}<s_{1}<\frac{3}{4}, \ s_{1}<s_{2}<1}$.

Let us use a closed ball in the Sobolev space, namely
\begin{equation}
\label{b}
B_{\rho}:=\{u(x)\in H^{2}({\mathbb R}^{3}) \ | \ \|u\|_{H^{2}({\mathbb R}^{3})}\leq
\rho \}, \quad 0<\rho\leq 1.
\end{equation}
We look for the solution of equation (\ref{pert}) as the fixed point of the
auxiliary nonlinear problem
\begin{equation}
\label{aux}
[(-\Delta)^{s_{1}}+(-\Delta)^{s_{2}}]u(x)=\varepsilon \int_{{\mathbb R}^{3}}
{\cal K}(x-y)g(u_{0}(y)+v(y))dy
\end{equation}
with
$\displaystyle{\frac{1}{4}<s_{1}<\frac{3}{4}, \ s_{1}<s_{2}<1}$
in ball (\ref{b}). For a given function $v(y)$ this is an equation with
respect to $u(x)$.

The left side of (\ref{aux}) involves the operator, which fails to satisfy
the Fredholm property
\begin{equation}
\label{l}
l:=(-\Delta)^{s_{1}}+(-\Delta)^{s_{2}}:
H^{2s_{2}}({\mathbb R}^{3})\to L^{2}({\mathbb R}^{3}).
\end{equation}
It is defined via the spectral calculus. This is the pseudo-differential
operator with the symbol $|p|^{2s_{1}}+|p|^{2s_{2}}$, so that
$$
lu(x)=\frac{1}{(2\pi)^{\frac{3}{2}}}\int_{{\mathbb R}^{3}}(|p|^{2s_{1}}+|p|^{2s_{2}})
\widehat{u}(p)e^{ipx}dp, \quad u(x)\in H^{2s_{2}}({\mathbb R}^{3}),
$$
with the standard Fourier transform is introduced in (\ref{f}).

Clearly, the essential spectrum of (\ref{l}) fills the nonnegative semi-axis
$[0, +\infty)$.
Therefore, this operator has no bounded inverse. The similar situation
appeared in works ~\cite{VV111} and ~\cite{VV14}. But as distinct from the
present situation, the solvability of the equations treated there required the
orthogonality conditions.

The fixed point technique was applied in ~\cite{V2010} to estimate the
perturbation to the standing solitary wave of the Nonlinear Schr\"odinger
(NLS) equation when either the external potential or the nonlinear term in the
NLS were perturbed but the Schr\"odinger operator involved in the nonlinear
problem there had the Fredholm property (see Assumption 1 of ~\cite{V2010},
also ~\cite{CPV05}).

Let us introduce the interval on the real line
\begin{equation}
\label{i}
I:=\big[-c_{e}\|u_{0}\|_{H^{2}({\mathbb R}^{3})}-c_{e}, \
c_{e}\|u_{0}\|_{H^{2}({\mathbb R}^{3})}+c_{e}\big]
\end{equation}
along with the closed ball in the space of $C^{2}(I)$ functions, namely
\begin{equation}
\label{M}
D_{M}:=\{g(z)\in C^{2}(I) \ | \ \|g\|_{C^{2}(I)}\leq M \}, \quad M>0.
\end{equation}
The norm contained in (\ref{M}) is given by
\begin{equation}
\label{gn}
\|g\|_{C^{2}(I)}:=\|g\|_{C(I)}+\|g'\|_{C(I)}+\|g''\|_{C(I)}
\end{equation}
with $\|g\|_{C(I)}:=\hbox{max}_{z\in I}|g(z)|$.

We impose the following technical conditions on the nonlinear part of
equation (\ref{p}). It will be trivial at the origin along with its first
derivative. From the perspective of the biological applications, $g(z)$ can
be, for example the quadratic function describing the cell-cell interaction.

\bigskip

\noindent
{\bf Assumption 1.2.} {\it Suppose that
$g(z): {\mathbb R}\to {\mathbb R}$, such that
$g(0)=0$ and $g'(0)=0$. Let us also assume that $g(z)\in D_{M}$ and
it does not vanish identically on the interval $I$}.

\bigskip

We introduce the operator $t_g$, so that $u = t_g v$, where $u$ is a
solution of equation (\ref{aux}). Our first main statement is as follows.

\bigskip

\noindent
{\bf Theorem 1.3.} {\it Suppose that Assumptions 1.1 and 1.2 are valid. Then
for every $\rho\in (0, 1]$ problem (\ref{aux})
defines the map $t_{g}: B_{\rho}\to B_{\rho}$, which is a strict contraction
for all
\begin{equation}
\label{eps}
0<\varepsilon   \leq \frac{\rho}{H} \end{equation}
where
$$ H = 2M(\|u_{0}\|_{H^{2}({\mathbb R}^{3})}+1)^{2}
\Big[\frac{\|{\cal K}\|_{L^{1}({\mathbb R}^{3})}^{2}
(\|u_{0}\|_{H^{2}({\mathbb R}^{3})}+1)^{\frac{8s_{1}}{3}-2}3}
{(3-4s_{1}){s_{1}}^{\frac{4s_{1}}{3}}{\pi}^{\frac{8s_{1}}{3}}{2}^{4s_{1}+\frac{8s_{1}}{3}}}
+\frac{Q^{2}}{4}\Big]^{\frac{1}{2}} . $$

The unique fixed point $u_{p}(x)$ of this map $t_{g}$ is the only solution of
equation (\ref{pert}) in $B_{\rho}$.}

\bigskip

\noindent
Note that the resulting solution of problem (\ref{p}) given by formula
(\ref{r}) will not vanish identically in the whole space because the source
term $f(x)$ is nontrivial and
$g(0)=0$ as we assume.

\noindent
We have the following elementary proposition.

\bigskip

\noindent
{\bf Lemma 1.4.} {\it For $R\in (0, +\infty)$  consider the
function
$$
\varphi(R):=\alpha R^{3-4s_{1}}+\frac{1}{R^{4s_{1}}}, \quad
\frac{1}{4}<s_{1}<\frac{3}{4}, \quad \alpha>0.
$$
It attains the minimal value at \
$\displaystyle{R^{*}:=\Bigg(\frac{4s_{1}}{\alpha (3-4s_{1})}\Bigg)^{\frac{1}{3}}}$,
which is given by}
$$
\varphi(R^{*})=3(3-4s_{1})^{\frac{4s_{1}}{3}-1}(4s_{1})^{-\frac{4s_{1}}{3}}
{\alpha}^{\frac{4s_{1}}{3}}.
$$

\bigskip

\noindent
Our second main statement deals with the continuity of the resulting solution
of problem (\ref{p}) given by (\ref{r}) with respect to
the nonlinear function $g$.

\bigskip

\noindent
{\bf Theorem 1.5.} {\it Suppose that $j=1,2$, the conditions of Theorem 1.3
are fulfilled, such that $u_{p,j}(x)$ is the unique fixed point of the map
$t_{g_{j}}: B_{\rho}\to B_{\rho}$, which is a strict contraction for all the values
of $\varepsilon$ satisfying (\ref{eps})
and the resulting solution of equation (\ref{p}) with
$g(z)=g_{j}(z)$ is given by
\begin{equation}
\label{ressol}
u_{j}(x)=u_{0}(x)+u_{p,j}(x).
\end{equation}
Then for all the values of $\varepsilon$, which satisfy bound (\ref{eps}),
the estimate
$$
\|u_{1}-u_{2}\|_{H^{2}({\mathbb R}^{3})}\leq \frac{\varepsilon}{1-\varepsilon \sigma}
(\|u_{0}\|_{H^{2}({\mathbb R}^{3})}+1)^{2}\times
$$
\begin{equation}
\label{cont}
\times \Bigg[
\frac{\|{\cal K}\|_{L^{1}({\mathbb R}^{3})}^{2}(\|u_{0}\|_{H^{2}
({\mathbb R}^{3})}+1)^{\frac{8s_{1}}{3}-2}}
{{s_{1}}^{\frac{4s_{1}}{3}}{\pi}^{\frac{8s_{1}}{3}}{2}^{4s_{1}+\frac{8s_{1}}{3}}}
\frac{3}{3-4s_{1}}+\frac{Q^{2}}{4}\Bigg]^
\frac{1}{2}\|g_{1}-g_{2}\|_{C^{2}(I)}
\end{equation}
holds.}

\bigskip

\noindent
Inequality (\ref{cont}) contains
$\sigma$, which is introduced in formula (\ref{sig}) further down.

\noindent
Let us proceed to the proof of our first main proposition.

\bigskip


\setcounter{equation}{0}

\section{\bf The existence of the perturbed solution}

\bigskip

\noindent
{\it Proof of Theorem 1.3.} Let us choose arbitrarily $v(x)\in B_{\rho}$ and
designate the term involved in the integral expression in the right side of
problem (\ref{aux}) as
$$
G(x):=g(u_{0}(x)+v(x)).
$$
We use the standard Fourier transform
\begin{equation}
\label{f}
\widehat{\phi}(p):=\frac{1}{(2\pi)^{\frac{3}{2}}}\int_{{\mathbb R}^{3}}\phi(x)
e^{-ipx}dx.
\end{equation}
Obviously, the estimate from above
\begin{equation}
\label{fub}
\|\widehat{\phi}(p)\|_{L^{\infty}({\mathbb R}^{3})}\leq \frac{1}{(2\pi)^{\frac{3}{2}}}
\|\phi(x)\|_{L^{1}({\mathbb R}^{3})}
\end{equation}
is valid. We apply (\ref{f}) to both sides of equation (\ref{aux}) and
arrive at
$$
\widehat{u}(p)=\varepsilon (2\pi)^{\frac{3}{2}}
\frac{\widehat{\cal K}(p)\widehat{G}(p)}{|p|^{2s_{1}}+|p|^{2s_{2}}},
$$
where
$\displaystyle{\frac{1}{4}<s_{1}<\frac{3}{4}, \ s_{1}<s_{2}<1}$.
Hence, for the norm we obtain
\begin{equation}
\label{un}
\|u\|_{L^{2}({\mathbb R}^{3})}^{2}={(2\pi)}^{3} \varepsilon^{2}\int_{{\mathbb R}^{3}}
\frac{|\widehat{\cal K}(p)|^{2}|\widehat{G}(p)|^{2}}
{[|p|^{2s_{1}}+|p|^{2s_{2}}]^{2}}dp.
\end{equation}
Note that as distinct from the previous works  ~\cite{VV111} and ~\cite{VV14}
having the
standard Laplacian in the diffusion term, here we do not try to control
the norm
$$
\Bigg\|\frac{\widehat{\cal K}(p)}{|p|^{2s_{1}}+|p|^{2s_{2}}}\Bigg\|_
{L^{\infty}({\mathbb R}^{3})}.
$$
Instead, we derive the bound on the right side of (\ref{un}) by means of the
analog of inequality (\ref{fub}) applied to functions ${\cal K}$ and $G$ with
$R>0$ as
$$
{(2\pi)}^{3} \varepsilon^{2}\int_{{\mathbb R}^{3}}
\frac{|\widehat{\cal K}(p)|^{2}|\widehat{G}(p)|^{2}}
{[|p|^{2s_{1}}+|p|^{2s_{2}}]^{2}}dp\leq
$$
$$
\leq {(2\pi)}^{3} \varepsilon^{2}\int_{|p|\leq R}\frac{|\widehat{\cal K}(p)|^{2}|
\widehat{G}(p)|^{2}}{|p|^{4s_{1}}}dp+
{(2\pi)}^{3} \varepsilon^{2}\int_{|p|>R}\frac{|\widehat{\cal K}(p)|^{2}
|\widehat{G}(p)|^{2}}{|p|^{4s_{1}}}dp\leq
$$
\begin{equation}
\label{ub1}
\leq \varepsilon^{2}\|{\cal K}\|_{L^{1}({\mathbb R}^{3})}^{2}\Bigg\{
\frac{1}{2{\pi}^{2}}\|G(x)\|_{L^{1}({\mathbb R}^{3})}^{2}\frac{R^{3-4s_{1}}}
{3-4s_{1}}+
\frac{1}{R^{4s_{1}}}\|G(x)\|_{L^{2}({\mathbb R}^{3})}^{2}\Bigg\}.
\end{equation}
For $v(x)\in B_{\rho}$, the estimate
$$
\|u_{0}+v\|_{L^{2}({\mathbb R}^{3})}\leq \|u_{0}\|_{H^{2}({\mathbb R}^{3})}+1
$$
is valid.
Sobolev embedding (\ref{e}) yields
$$
|u_{0}+v|\leq c_{e}(\|u_{0}\|_{H^{2}({\mathbb R}^{3})}+1).
$$
Evidently,
$$
G(x)=\int_{0}^{u_{0}+v}g'(z)dz,
$$
so that
$$
|G(x)|\leq \hbox{sup}_{z\in I}|g'(z)||u_{0}+v|\leq M|u_{0}+v|
$$
with the interval $I$ defined in (\ref{i}).
Therefore,
$$
\|G(x)\|_{L^{2}({\mathbb R}^{3})}\leq M\|u_{0}+v\|_{L^{2}({\mathbb R}^{3})}\leq M
(\|u_{0}\|_{H^{2}({\mathbb R}^{3})}+1).
$$
Clearly,
$$
G(x)=\int_{0}^{u_{0}+v}dy\Big[\int_{0}^{y}g''(z)dz \Big].
$$
This implies
$$
|G(x)|\leq \frac{1}{2}\hbox{sup}_{z\in I}|g''(z)||u_{0}+v|^{2}\leq
\frac{M}{2}|u_{0}+v|^{2},
$$
such that
\begin{equation}
\label{G1}
\|G(x)\|_{L^{1}({\mathbb R}^{3})}\leq \frac{M}{2}\|u_{0}+v\|_{L^{2}({\mathbb R}^{3})}^{2}
\leq \frac{M}{2}(\|u_{0}\|_{H^{2}({\mathbb R}^{3})}+1)^{2}.
\end{equation}
Thus, we obtain the upper bound for the right side of (\ref{ub1}) equal to
$$
{\varepsilon}^{2}\|{\cal K}\|_{L^{1}({\mathbb R}^{3})}^{2}M^{2}
(\|u_{0}\|_{H^{2}({\mathbb R}^{3})}+1)^{2}
\Bigg\{\frac{(\|u_{0}\|_{H^{2}({\mathbb R}^{3})}+1)^{2}R^{3-4s_{1}}}{8\pi^{2}
(3-4s_{1})}+\frac{1}{R^{4s_{1}}}\Bigg\},
$$
where $R\in (0, +\infty)$. Let us use Lemma 1.4 to derive the minimal value
of such expression. Hence, $\|u\|_{L^{2}({\mathbb R}^{3})}^{2}\leq$
\begin{equation}
\label{ul2ub}
\leq {\varepsilon}^{2}\|{\cal K}\|
_{L^{1}({\mathbb R}^{3})}^{2}M^{2}(\|u_{0}\|_{H^{2}({\mathbb R}^{3})}+1)^{2+\frac{8s_{1}}{3}}
\frac{3}{(3-4s_{1}){s_{1}}^{\frac{4s_{1}}{3}}{\pi}^{\frac{8s_{1}}{3}}
2^{4s_{1}+\frac{8s_{1}}{3}}}.
\end{equation}
Obviously, by means of (\ref{aux}) we have
$$
[-\Delta+(-\Delta)^{1+s_{2}-s_{1}}]u(x)=\varepsilon
(-\Delta)^{1-s_{1}}\int_{{\mathbb R}^{3}}{\cal K}(x-y)G(y)dy.
$$
By virtue of the standard Fourier transform (\ref{f}), the analog of bound
(\ref{fub}) applied to function $G$ and inequality (\ref{G1}) we arrive at
\begin{equation}
\label{32}
\|\Delta u\|_{L^{2}({\mathbb R}^{3})}^{2}\leq \varepsilon^{2}\|G\|_
{L^{1}({\mathbb R}^{3})}^{2}Q^{2}\leq \varepsilon^{2}\frac{M^{2}}{4}
(\|u_{0}\|_{H^{2}({\mathbb R}^{3})}+1)^{4}Q^{2}.
\end{equation}
Let us recall the definition of the norm (\ref{h3}).
Estimates (\ref{ul2ub}) and (\ref{32}) give us
$$
\|u\|_{H^{2}({\mathbb R}^{3})}\leq \varepsilon (\|u_{0}\|_{H^{2}({\mathbb R}^{3})}+1)^{2}M
\times
$$
\begin{equation}
\label{uh3e}
\times
 \Bigg[\frac{\|{\cal K}\|_{L^{1}({\mathbb R}^{3})}^{2}(\|u_{0}\|_{H^{2}
({\mathbb R}^{3})}+1)^{\frac{8s_{1}}{3}-2}3}{(3-4s_{1})s_{1}^{\frac{4s_{1}}{3}}
{\pi}^{\frac{8s_{1}}{3}}{2}^{4s_{1}+\frac{8s_{1}}{3}}}
+\frac{Q^{2}}{4}\Bigg]^{\frac{1}{2}}\leq \rho
\end{equation}
for all the values of the parameter $\varepsilon$, which satisfy bound
(\ref{eps}). This means that $u(x)\in B_{\rho}$ as well.

Suppose that for some $v(x)\in B_{\rho}$ there exist
two solutions $u_{1,2}(x)\in B_{\rho}$ of problem (\ref{aux}). Evidently, the
difference function $w(x):=u_{1}(x)-u_{2}(x) \in H^{2}({\mathbb R}^{3})$ is a
solution of the homogeneous equation
$$
[(-\Delta)^{s_{1}}+(-\Delta)^{s_{2}}]w=0.
$$
Because the operator $l: H^{2s_{2}}({\mathbb R}^{3})\to L^{2}({\mathbb R}^{3})$
given by (\ref{l}) does not have any
nontrivial zero modes, $w(x)$ will vanish in the whole space.
Therefore, problem (\ref{aux}) defines a map
$t_{g}: B_{\rho}\to B_{\rho}$ for all the values of $\varepsilon$
satisfying (\ref{eps}).

\noindent
Our goal is to demonstrate that under the given conditions such map is a strict
contraction. Let us choose arbitrary $v_{1,2}(x)\in B_{\rho}$. The reasoning above
yields
$u_{1,2}:=t_{g}v_{1,2}\in B_{\rho}$ as well if $\varepsilon$ satisfies
inequality (\ref{eps}). Using (\ref{aux}), we have
\begin{equation}
\label{aux1}
[(-\Delta)^{s_{1}}+(-\Delta)^{s_{2}}]u_{1}(x)=\varepsilon \int_{{\mathbb R}^{3}}
{\cal K}(x-y)g(u_{0}(y)+v_{1}(y))dy,
\end{equation}
\begin{equation}
\label{aux2}
[(-\Delta)^{s_{1}}+(-\Delta)^{s_{2}}]u_{2}(x)=\varepsilon \int_{{\mathbb R}^{3}}
{\cal K}(x-y)g(u_{0}(y)+v_{2}(y))dy
\end{equation}
with $\displaystyle{\frac{1}{4}<s_{1}<\frac{3}{4}, \ s_{1}<s_{2}<1}$.
We define
$$
G_{1}(x):=g(u_{0}(x)+v_{1}(x)), \quad G_{2}(x):=g(u_{0}(x)+v_{2}(x)).
$$
Let us apply the standard Fourier transform (\ref{f}) to both sides of
equations (\ref{aux1}) and (\ref{aux2}) to obtain
$$
\widehat{u_{1}}(p)=\varepsilon (2\pi)^{\frac{3}{2}}
\frac{\widehat{\cal K}(p)\widehat{G_{1}}(p)}{|p|^{2s_{1}}+|p|^{2s_{2}}}, \quad
\widehat{u_{2}}(p)=\varepsilon (2\pi)^{\frac{3}{2}}
\frac{\widehat{\cal K}(p)\widehat{G_{2}}(p)}{|p|^{2s_{1}}+|p|^{2s_{2}}},
$$
so that
\begin{equation}
\label{u12n2d}
\|u_{1}-u_{2}\|_{L^{2}({\mathbb R}^{3})}^{2}=\varepsilon^{2}{(2\pi)}^{3}
\int_{{\mathbb R}^{3}}\frac{|\widehat{\cal K}(p)|^{2}
|{\widehat{G_{1}}(p)}-{\widehat{G_{2}}(p)}|^{2}}{[|p|^{2s_{1}}+|p|^{2s_{2}}]^{2}}dp.
\end{equation}
The right side of (\ref{u12n2d}) can be trivially estimated from above by
means of bound (\ref{fub}) as
$$
\varepsilon^{2}{(2\pi)}^{3}\int_{|p|\leq R}\frac{|\widehat{\cal K}(p)|^{2}
|{\widehat{G_{1}}(p)}-{\widehat{G_{2}}(p)}|^{2}}{|p|^{4s_{1}}}dp+
$$
$$
+\varepsilon^{2}{(2\pi)}^{3}\int_{|p|>R}\frac{|\widehat{\cal K}(p)|^{2}
|{\widehat{G_{1}}(p)}-{\widehat{G_{2}}(p)}|^{2}}{|p|^{4s_{1}}}dp \leq
$$
$$
\leq \varepsilon^{2}\|{\cal K}\|_{L^{1}({\mathbb R}^{3})}^{2}\Bigg\{\frac
{\|G_{1}(x)-G_{2}(x)\|_{L^{1}({\mathbb R}^{3})}^{2}}{2\pi^{2}}\frac{R^{3-4s_{1}}}{3-4s_{1}}+
\frac{\|G_{1}(x)-G_{2}(x)\|_{L^{2}({\mathbb R}^{3})}^{2}}{R^{4s_{1}}}\Bigg\}
$$
with $R\in (0,+\infty)$. Note that
$$
G_{1}(x)-G_{2}(x)=\int_{u_{0}+v_{2}}^{u_{0}+v_{1}}g'(z)dz.
$$
Thus,
$$
|G_{1}(x)-G_{2}(x)|\leq \hbox{sup}_{z\in I}|g'(z)||v_{1}(x)-v_{2}(x)|\leq
M|v_{1}(x)-v_{2}(x)|,
$$
such that
$$
\|G_{1}(x)-G_{2}(x)\|_{L^{2}({\mathbb R}^{3})}\leq M\|v_{1}-v_{2}\|_
{L^{2}({\mathbb R}^{3})}\leq M\|v_{1}-v_{2}\|_{H^{2}({\mathbb R}^{3})}.
$$
Let us write
$$
G_{1}(x)-G_{2}(x)=\int_{u_{0}+v_{2}}^{u_{0}+v_{1}}dy \Big[\int_{0}^{y}g''(z)dz \Big].
$$
Clearly, $G_{1}(x)-G_{2}(x)$ can be easily bounded from above in the absolute
value by
$$
\frac{1}{2}\hbox{sup}_{z\in I}|g''(z)||(v_{1}-v_{2})(2u_{0}+
v_{1}+v_{2})|\leq\frac{M}{2}|(v_{1}-v_{2})(2u_{0}+v_{1}+v_{2})|.
$$
We use the Schwarz inequality to derive the estimate from above for the norm
$$
\|G_{1}(x)-G_{2}(x)\|_{L^{1}({\mathbb R}^{3})}\leq
\frac{M}{2}\|v_{1}-v_{2}\|_{L^{2}({\mathbb R}^{3})}\|2u_{0}+v_{1}+v_{2}\|_
{L^{2}({\mathbb R}^{3})}\leq
$$
\begin{equation}
\label{g12}
\leq
M\|v_{1}-v_{2}\|_{H^{2}({\mathbb R}^{3})}
(\|u_{0}\|_{H^{2}({\mathbb R}^{3})}+1).
\end{equation}
Hence, we obtain the upper bound for
$\|u_{1}(x)-u_{2}(x)\|_{L^{2}({\mathbb R}^{3})}^{2}$ equal to
$$
\varepsilon^{2}\|{\cal K}\|
_{L^{1}({\mathbb R}^{3})}^{2}M^{2}\|v_{1}-v_{2}\|_{H^{2}({\mathbb R}^{3})}^{2}
\Big\{\frac{(\|u_{0}\|_{H^{2}({\mathbb R}^{3})}+1)^{2}}{2{\pi}^{2}}
\frac{R^{3-4s_{1}}}{3-4s_{1}}+\frac{1}{R^{4s_{1}}}\Big\}.
$$
Let us recall Lemma 1.4 and minimize such expression over
$R\in (0,+\infty)$. Thus,
$$
\|u_{1}(x)-u_{2}(x)\|_{L^{2}({\mathbb R}^{3})}^{2}\leq
$$
\begin{equation}
\label{u12n}
\varepsilon^{2}\|{\cal K}\|_{L^{1}({\mathbb R}^{3})}^{2}M^{2}
\|v_{1}-v_{2}\|_{H^{2}({\mathbb R}^{3})}^{2}\frac{3}{3-4s_{1}}
\frac{(\|u_{0}\|_{H^{2}({\mathbb R}^{3})}+1)^{\frac{8s_{1}}{3}}}
{{s_{1}}^{\frac{4s_{1}}{3}}{\pi}^{\frac{8s_{1}}{3}}2^{4s_{1}}}.
\end{equation}
Formulas (\ref{aux1}) and (\ref{aux2}) give us
$$
[-\Delta+(-\Delta)^{1+s_{2}-s_{1}}](u_{1}(x)-u_{2}(x))=\varepsilon
(-\Delta)^{1-s_{1}}
\int_{{\mathbb R}^{3}}{\cal K}(x-y)[G_{1}(y)-G_{2}(y)]dy.
$$
By means of the standard Fourier transform (\ref{f}) along with estimates
(\ref{fub}) and (\ref{g12}), we arrive at
$$
\|\Delta(u_{1}-u_{2})\|_{L^{2}({\mathbb R}^{3})}^{2}\leq \varepsilon^{2}
Q^{2}\|G_{1}-G_{2}\|_{L^{1}({\mathbb R}^{3})}^{2}\leq
$$
\begin{equation}
\label{d12}
\leq \varepsilon^{2}Q^{2}M^{2}
\|v_{1}-v_{2}\|_{H^{2}({\mathbb R}^{3})}^{2}(\|u_{0}\|_{H^{2}({\mathbb R}^{3})}+1)^{2}.
\end{equation}
By virtue of inequalities (\ref{u12n}) and (\ref{d12}), the norm
$\|u_{1}-u_{2}\|_{H^{2}({\mathbb R}^{3})}$ can be bounded from above by the quantity
$$
\varepsilon M(\|u_{0}\|_{H^{2}({\mathbb R}^{3})}+1)\times
$$
$$
\times
\Bigg\{\frac{\|{\cal K}\|_{L^{1}({\mathbb R}^{3})}^{2}
(\|u_{0}\|_{H^{2}({\mathbb R}^{3})}+1)^{\frac{8s_{1}}{3}-2}}
{{s_{1}}^{\frac{4s_{1}}{3}}{\pi}^{\frac{8s_{1}}{3}}2^{4s_{1}}}\frac{3}{3-4s_{1}}+Q^{2}\Bigg\}^{\frac{1}{2}}
\times
$$
\begin{equation}
\label{contr}
\times\|v_{1}-v_{2}\|_{H^{2}({\mathbb R}^{3})}.
\end{equation}
It can be trivially checked using (\ref{eps}) that the constant in the right
side of (\ref{contr}) is less than one.
Therefore, the map $t_{g}: B_{\rho}\to B_{\rho}$ defined by problem
(\ref{aux}) is a strict contraction for all the values of $\varepsilon$,
which satisfy bound (\ref{eps}). Its unique fixed point $u_{p}(x)$ is the
only solution of equation (\ref{pert}) in the ball $B_{\rho}$. By means of
(\ref{uh3e}), we obtain that $\|u_{p}(x)\|_{H^{2}({\mathbb R}^{3})}\to 0$ as
$\varepsilon\to 0$.
The resulting
$u(x)\in H^{2}({\mathbb R}^{3})$ given by (\ref{r}) is a solution of
problem (\ref{p}). \hfill\lanbox

\bigskip

We turn our attention to demonstrating the validity of the second main
statement of the work.

\bigskip


\setcounter{equation}{0}

\section{\bf The continuity of the resulting solution}

\bigskip

\noindent
{\it Proof of Theorem 1.5.} Obviously, for all the values of the parameter
$\varepsilon$ satisfying inequality (\ref{eps}),  we have
$$
u_{p,1}=t_{g_{1}}u_{p,1}, \quad u_{p,2}=t_{g_{2}}u_{p,2},
$$
such that
$$
u_{p,1}-u_{p,2}=t_{g_{1}}u_{p,1}-t_{g_{1}}u_{p,2}+t_{g_{1}}u_{p,2}-
t_{g_{2}}u_{p,2}.
$$
Evidently,
$$
\|u_{p,1}-u_{p,2}\|_{H^{2}({\mathbb R}^{3})}\leq\|t_{g_{1}}u_{p,1}-t_{g_{1}}
u_{p,2}\|_{H^{2}({\mathbb R}^{3})}+\|t_{g_{1}}u_{p,2}-t_{g_{2}}u_{p,2}\|_
{H^{2}({\mathbb R}^{3})}.
$$
Let us define the positive constant
$$
\sigma:=M(\|u_{0}\|_{H^{2}({\mathbb R}^{3})}+1)\times
$$
\begin{equation}
\label{sig}
\times
\Bigg\{\frac{\|{\cal K}\|_{L^{1}({\mathbb R}^{3})}^{2}(\|u_{0}\|
_{H^{2}({\mathbb R}^{3})}+1)^{\frac{8s_{1}}{3}-2}}
{{s_{1}}^{\frac{4s_{1}}{3}}{\pi}^{\frac{8s_{1}}{3}}2^{4s_{1}}}\frac{3}{3-4s_{1}}+Q^{2}
\Bigg\}^{\frac{1}{2}}.
\end{equation}
By virtue of estimate from above (\ref{contr}), we derive
$$
\|t_{g_{1}}u_{p,1}-t_{g_{1}}u_{p,2}\|_{H^{2}({\mathbb R}^{3})}\leq \varepsilon
\sigma\|u_{p,1}-u_{p,2}\|_{H^{2}({\mathbb R}^{3})}.
$$
Clearly, $\varepsilon \sigma<1$ since the map
$t_{g_{1}}: B_{\rho}\to  B_{\rho}$ is a strict contraction under the stated
assumptions.
Thus,
\begin{equation}
\label{sigma}
(1-\varepsilon \sigma)\|u_{p,1}-u_{p,2}\|_{H^{2}({\mathbb R}^{3})}\leq
\|t_{g_{1}}u_{p,2}-t_{g_{2}}u_{p,2}\|_{H^{2}({\mathbb R}^{3})}.
\end{equation}
Let us recall that for our fixed point $t_{g_{2}}u_{p,2}=u_{p,2}$ and define
$\gamma(x):=t_{g_{1}}u_{p,2}$. Obviously,
\begin{equation}
\label{12}
[(-\Delta)^{s_{1}}+(-\Delta)^{s_{2}}]\gamma(x)=\varepsilon \int_{{\mathbb R}^{3}}
{\cal K}(x-y)g_{1}(u_{0}(y)+u_{p,2}(y))dy,
\end{equation}
\begin{equation}
\label{22}
[(-\Delta)^{s_{1}}+(-\Delta)^{s_{2}}]u_{p,2}(x)=\varepsilon \int_{{\mathbb R}^{3}}
{\cal K}(x-y)g_{2}(u_{0}(y)+u_{p,2}(y))dy
\end{equation}
with
$\displaystyle{\frac{1}{4}<s_{1}<\frac{3}{4}, \ s_{1}<s_{2}<1}$.
We introduce
$$
G_{1, 2}(x):=g_{1}(u_{0}(x)+u_{p,2}(x)), \quad G_{2, 2}(x):=g_{2}(u_{0}(x)+u_{p,2}(x)).
$$
Let us apply the standard Fourier transform (\ref{f}) to both sides of
equations (\ref{12}) and (\ref{22}). This gives us
$$
\widehat{\gamma}(p)=\varepsilon (2 \pi)^{\frac{3}{2}} \frac{\widehat{\cal K}(p)
\widehat{G_{1, 2}}(p)}{|p|^{2s_{1}}+|p|^{2s_{2}}}, \quad
\widehat{u_{p,2}}(p)=\varepsilon (2 \pi)^{\frac{3}{2}} \frac{\widehat{\cal K}(p)
\widehat{G_{2, 2}}(p)}{|p|^{2s_{1}}+|p|^{2s_{2}}}.
$$
Evidently,
\begin{equation}
\label{xiup2n}
\|\gamma(x)-u_{p,2}(x)\|_{L^{2}({\mathbb R}^{3})}^{2}=\varepsilon^{2}{(2\pi)}^{3}
\int_{{\mathbb R}^{3}}\frac{|{\widehat{\cal K}}(p)|^{2}|\widehat{G_{1, 2}}(p)-
\widehat{G_{2, 2}}(p)|^{2}}{[|p|^{2s_{1}}+|p|^{2s_{2}}]^{2}}dp.
\end{equation}
We use (\ref{fub}) to derive the upper bound on
the right side of (\ref{xiup2n}), namely
$$
\varepsilon^{2}{(2\pi)}^{3}
\int_{|p|\leq R}\frac{|{\widehat{\cal K}}(p)|^{2}|\widehat{G_{1, 2}}(p)-
\widehat{G_{2, 2}}(p)|^{2}}{|p|^{4s_{1}}}dp+
$$
$$
+\varepsilon^{2}{(2\pi)}^{3}
\int_{|p|>R}\frac{|{\widehat{\cal K}}(p)|^{2}|\widehat{G_{1, 2}}(p)-
\widehat{G_{2, 2}}(p)|^{2}}{|p|^{4s_{1}}}dp\leq
$$
$$
\leq \varepsilon^{2}\|{\cal K}\|_{L^{1}({\mathbb R}^{3})}^{2}\Bigg\{\frac{1}
{2{\pi}^{2}}\|G_{1, 2}-G_{2, 2}\|_{L^{1}({\mathbb R}^{3})}^{2}\frac{R^{3-4s_{1}}}{3-4s_{1}}+
\frac{\|G_{1, 2}-G_{2, 2}\|_{L^{2}({\mathbb R}^{3})}^{2}}{R^{4s_{1}}}\Bigg\},
$$
where $R\in (0, +\infty)$. We have
$$
G_{1, 2}(x)-G_{2, 2}(x)=\int_{0}^{u_{0}(x)+u_{p,2}(x)}[g_{1}'(z)-g_{2}'(z)]dz,
$$
so that
$$
|G_{1, 2}(x)-G_{2, 2}(x)|\leq \hbox{sup}_{z\in I}|g_{1}'(z)-g_{2}'(z)|
|u_{0}(x)+u_{p,2}(x)|\leq
$$
$$
\leq\|g_{1}-g_{2}\|_{C^{2}(I)}|u_{0}(x)+u_{p,2}(x)|.
$$
Hence,
$$
\|G_{1, 2}-G_{2, 2}\|_{L^{2}({\mathbb R}^{3})}\leq \|g_{1}-g_{2}\|_{C^{2}(I)}
\|u_{0}+u_{p,2}\|_{L^{2}({\mathbb R}^{3})}\leq
$$
$$
\leq \|g_{1}-g_{2}\|_{C^{2}(I)}(\|u_{0}\|_{H^{2}({\mathbb R}^{3})}+1).
$$
We make use of another equality
$$
G_{1, 2}(x)-G_{2, 2}(x)=\int_{0}^{u_{0}(x)+u_{p,2}(x)}dy\Big[\int_{0}^{y}
(g_{1}''(z)-g_{2}''(z))dz\Big].
$$
Clearly,
$$
|G_{1, 2}(x)-G_{2, 2}(x)|\leq \frac{1}{2}\hbox{sup}_{z\in I}|g_{1}''(z)-g_{2}''(z)|
|u_{0}(x)+u_{p,2}(x)|^{2}\leq
$$
$$
\leq  \frac{1}{2}\|g_{1}-g_{2}\|_{C^{2}(I)}|u_{0}(x)+u_{p,2}(x)|^{2}.
$$
Therefore,
$$
\|G_{1, 2}-G_{2, 2}\|_{L^{1}({\mathbb R}^{3})}\leq \frac{1}{2}
\|g_{1}-g_{2}\|_{C^{2}(I)}\|u_{0}+u_{p,2}\|_{L^{2}({\mathbb R}^{3})}^{2}\leq
$$
\begin{equation}
\label{G1222}
\leq \frac{1}{2}\|g_{1}-g_{2}\|_{C^{2}(I)}(\|u_{0}\|_{H^{2}({\mathbb R}^{3})}+1)^{2}.
\end{equation}
This enables us to estimate from above the norm
$\|\gamma(x)-u_{p,2}(x)\|_{L^{2}({\mathbb R}^{3})}^{2}$ by
$$
\varepsilon^{2}\|{\cal K}\|_{L^{1}({\mathbb R}^{3})}^{2}
(\|u_{0}\|_{H^{2}({\mathbb R}^{3})}+1)^{2}\|g_{1}-g_{2}\|_{C^{2}(I)}^{2}\times
$$
\begin{equation}
\label{xiup2n2}
\times\Big[\frac{1}{8{\pi}^{2}}(\|u_{0}\|_{H^{2}({\mathbb R}^{3})}+1)^{2}
\frac{R^{3-4s_{1}}}{3-4s_{1}}+\frac{1}{R^{4s_{1}}}\Big].
\end{equation}
We minimize expression (\ref{xiup2n2}) over $R\in (0, +\infty)$
via Lemma 1.4. Thus,
$$
\|\gamma(x)-u_{p,2}(x)\|_{L^{2}({\mathbb R}^{3})}^{2}\leq
$$
$$
\leq \varepsilon^{2}
\|{\cal K}\|_{L^{1}({\mathbb R}^{3})}^{2}
\|g_{1}-g_{2}\|_{C^{2}(I)}^{2}\frac{3}{3-4s_{1}}
\frac{(\|u_{0}\|_{H^{2}({\mathbb R}^{3})}+1)^{2+\frac{8s_{1}}{3}}}{{s_{1}}^{\frac{4s_{1}}{3}}
{\pi}^{\frac{8s_{1}}{3}}2^{4s_{1}+\frac{8s_{1}}{3}}}.
$$
Let us use formulas (\ref{12}) and (\ref{22}) to derive
$$
[-\Delta+(-\Delta)^{1+s_{2}-s_{1}}]\gamma(x)=\varepsilon
(-\Delta)^{1-s_{1}}\int_{{\mathbb R}^{3}}
{\cal K}(x-y)G_{1, 2}(y)dy,
$$
$$
[-\Delta+(-\Delta)^{1+s_{2}-s_{1}}]u_{p,2}(x)=\varepsilon
(-\Delta)^{1-s_{1}}\int_{{\mathbb R}^{3}}
{\cal K}(x-y)G_{2, 2}(y)dy.
$$
This gives us
$$
[-\Delta+(-\Delta)^{1+s_{2}-s_{1}}](\gamma(x)-u_{p,2}(x))=
$$
$$
=\varepsilon (-\Delta)^{1-s_{1}}\int_{{\mathbb R}^{3}}{\cal K}(x-y)
[G_{1, 2}(y)-G_{2, 2}(y)]dy.
$$
By means of the standard Fourier transform (\ref{f}) along with bounds
(\ref{fub}) and (\ref{G1222}), we obtain that
$$
\|\Delta [\gamma(x)-u_{p,2}(x)]\|_{L^{2}({\mathbb R}^{3})}^{2}\leq
$$
$$
\leq \varepsilon^{2}\|G_{1,2}-G_{2,2}\|_{L^{1}({\mathbb R}^{3})}^{2}Q^{2}\leq
\frac{\varepsilon^{2}Q^{2}}{4}(\|u_{0}\|_{H^{2}({\mathbb R}^{3})}+1)^{4}
\|g_{1}-g_{2}\|_{C^{2}(I)}^{2}.
$$
Thus,
$$
\|\gamma(x)-u_{p,2}(x)\|_{H^{2}({\mathbb R}^{3})}\leq \varepsilon
\|g_{1}-g_{2}\|_{C^{2}(I)}(\|u_{0}\|_{H^{2}({\mathbb R}^{3})}+1)^{2}\times
$$
$$
\times \Bigg[
\frac{\|{\cal K}\|_{L^{1}({\mathbb R}^{3})}^{2}(\|u_{0}\|_{H^{2}
({\mathbb R}^{3})}+1)^{\frac{8s_{1}}{3}-2}}
{{s_{1}}^{\frac{4s_{1}}{3}}{\pi}^{\frac{8s_{1}}{3}}{2}^{4s_{1}+\frac{8s_{1}}{3}}}
\frac{3}{3-4s_{1}}+\frac{Q^{2}}{4}\Bigg]^
\frac{1}{2}.
$$
Let us recall inequality (\ref{sigma}). Hence,
$$
\|u_{p,1}-u_{p,2}\|_{H^{2}({\mathbb R}^{3})}\leq
\frac{\varepsilon}{1-\varepsilon \sigma}(\|u_{0}\|_{H^{2}({\mathbb R}^{3})}+1)^{2}\times
$$
$$
\times \Bigg[
\frac{\|{\cal K}\|_{L^{1}({\mathbb R}^{3})}^{2}(\|u_{0}\|_{H^{2}
({\mathbb R}^{3})}+1)^{\frac{8s_{1}}{3}-2}}
{{s_{1}}^{\frac{4s_{1}}{3}}{\pi}^{\frac{8s_{1}}{3}}{2}^{4s_{1}+\frac{8s_{1}}{3}}}
\frac{3}{3-4s_{1}}+\frac{Q^{2}}{4}\Bigg]^
\frac{1}{2}\|g_{1}-g_{2}\|_{C^{2}(I)}.
$$
Formula (\ref{ressol}) allows us to complete the proof of the theorem.
\hfill\lanbox

\bigskip


\setcounter{equation}{0}

\section{Auxiliary results}

\bigskip

\noindent
We obtain the solvability relations for the linear
Poisson type problem with a square integrable right side in the situation
of the double scale anomalous diffusion
\begin{equation}
\label{lp}
[{(-\Delta)}^{s_{1}}+{(-\Delta)}^{s_{2}}]u=f(x), \quad x\in {\mathbb R}^{3},
\quad 0<s_{1}<s_{2}<1.
\end{equation}
The inner product is given by
\begin{equation}
\label{ip}
(f(x), g(x))_{L^{2} ({\mathbb R}^{3})}:=\int_{{\mathbb R}^{3}}f(x)\bar{g}(x)dx,
\end{equation}
with a slight abuse of notations when the functions contained in (\ref{ip})
are not square integrable, like for example the one involved in orthogonality
condition (\ref{oc1}) below. Indeed, if $f(x)\in L^{1} ({\mathbb R}^{3})$
and $g(x)\in L^{\infty}({\mathbb R}^{3})$, then the integral in the right side of
(\ref{ip}) is well defined.

The technical proposition below is easily established by
applying (\ref{f}) to both sides of equation
(\ref{lp}).

\bigskip

\noindent
{\bf Lemma 4.1.} {\it  Let  $0<s_{1}<s_{2}<1, \ f(x): {\mathbb R}^{3}\to
{\mathbb R}, \ f(x)\in L^{2} ({\mathbb R}^{3})$.

\medskip

\noindent
a) If $\displaystyle{s_{1}\in \Big(0, \frac{3}{4}\Big)}$ and additionally
$f(x)\in L^{1} ({\mathbb R}^{3})$, then problem (\ref{lp}) admits a unique
solution $u(x)\in H^{2s_{2}} ({\mathbb R}^{3})$.

\medskip

\noindent
b) If $\displaystyle{s_{1}\in \Big[\frac{3}{4}, 1\Big)}$ and in addition
$xf(x)\in L^{1} ({\mathbb R}^{3})$, then equation (\ref{lp}) has a unique
solution $u(x)\in H^{2s_{2}} ({\mathbb R}^{3})$ if and only if the orthogonality
relation
\begin{equation}
\label{oc1}
(f(x), 1)_{L^{2} ({\mathbb R}^{3})}=0
\end{equation}
holds.}

\bigskip

\noindent
{\it Proof.} It can be trivially checked that if
$u(x)\in L^{2} ({\mathbb R}^{3})$ solves problem (\ref{lp}) with
$f(x)\in L^{2}({\mathbb R}^{3})$, it will be contained in
$H^{2s_{2}}({\mathbb R}^{3})$ as well. To establish that, we apply the standard
Fourier transform (\ref{f}) to both sides of (\ref{lp}) and arrive at
$$
(|p|^{2s_{1}}+|p|^{2s_{2}})\widehat{u}(p)=\widehat{f}(p)\in L^{2} ({\mathbb R}^{3}).
$$
Hence,
$$
\int_{{\mathbb R}^{3}}[|p|^{2s_{1}}+|p|^{2s_{2}}]^{2}|\widehat{u}(p)|^{2}dp<\infty.
$$
Clearly, the trivial identity
$$
\|(-\Delta)^{s_{2}}u\|_{L^{2} ({\mathbb R}^{3})}^{2}=\int_{{\mathbb R}^{3}}|p|^{4s_{2}}
|\widehat{u}(p)|^{2}dp<\infty
$$
is valid.
Therefore,
$(-\Delta)^{s_{2}}u\in L^{2} ({\mathbb R}^{3})$. Using the definition
of the norm (\ref{n}), we obtain that
$u(x)\in H^{2s_{2}}({\mathbb R}^{3})$.

Let us address the uniqueness of solutions for our equation. Suppose that
problem (\ref{lp}) has two solutions
$u_{1, 2}(x)\in H^{2s_{2}}({\mathbb R}^{3})$. Then the difference function
$w(x):=u_{1}(x)-u_{2}(x)\in H^{2s_{2}}({\mathbb R}^{3})$ satisfies the homogeneous
equation
$$
[(-\Delta)^{s_{1}}+(-\Delta)^{s_{2}}]w=0.
$$
Since the operator
$l: H^{2s_{2}}({\mathbb R}^{3})\to L^{2}({\mathbb R}^{3})$ introduced in (\ref{l})
does not possess
any nontrivial zero modes, $w(x)$ will vanish in ${\mathbb R}^{3}$.

We apply the standard Fourier transform (\ref{f}) to both sides of problem
(\ref{lp}), which yields
\begin{equation}
\label{lpf}
\widehat{u}(p)=\frac{\widehat{f}(p)}{|p|^{2s_{1}}+|p|^{2s_{2}}}
\chi_{\{|p|\leq 1\}}+
\frac{\widehat{f}(p)}{|p|^{2s_{1}}+|p|^{2s_{2}}}
\chi_{\{|p|>1\}}.
\end{equation}
Here and further down $\chi_{A}$ will denote the characteristic function of a
set $A\subseteq {\mathbb R}^{3}$.

Evidently, the second term in the right side of
(\ref{lpf}) can be estimated from above in the absolute value by
$\displaystyle{\frac{|\widehat{f}(p)|}{2}\in L^{2}({\mathbb R}^{3})}$
according to the one of our assumptions.

Obviously, the first term in the right side of (\ref{lpf}) can be bounded
from above in the absolute value by
\begin{equation}
\label{fps2}
\frac{\|f(x)\|_{L^{1}({\mathbb R}^{3})}}{(2\pi)^{\frac{3}{2}}|p|^{2s_{1}}}\chi_{\{|p|\leq 1\}}
\end{equation}
using inequality (\ref{fub}). It can be easily verified that expression
(\ref{fps2}) belongs to $L^{2}({\mathbb R}^{3})$ for
$\displaystyle{s_{1}\in \Big(0, \frac{3}{4}\Big)}$.

\noindent
Let us turn our attention to the situation when
$\displaystyle{s_{1}\in \Big[\frac{3}{4}, 1\Big)}$. We write
\begin{equation}
\label{fhint}
\widehat{f}(p)=\widehat{f}(0)+\int_{0}^{|p|}
\frac{\partial \widehat{f}(q, \sigma)}{\partial q}dq.
\end{equation}
Here and below $\sigma$ will designate the angle variables on the sphere.
Thus, the first term in the right side of (\ref{lpf}) equals to
\begin{equation}
\label{fh0is12}
\frac{\widehat{f}(0)}{|p|^{2s_{1}}+|p|^{2s_{2}}}\chi_{\{|p|\leq 1\}}+
\frac{\int_{0}^{|p|}\frac{\partial \widehat{f}(q, \sigma)}{\partial q}dq}
{|p|^{2s_{1}}+|p|^{2s_{2}}}\chi_{\{|p|\leq 1\}}.
\end{equation}
By means of the definition of the standard Fourier transform (\ref{f}),
we easily obtain that
\begin{equation}
\label{pfhdp}
\Bigg|\frac{\partial \widehat{f}(p)}{\partial |p|}\Bigg|\leq
\frac{\|xf(x)\|_{L^{1}({\mathbb R}^{3})}}{(2\pi)^{\frac{3}{2}}}.
\end{equation}
This yields
$$
\Bigg|\frac{\int_{0}^{|p|}\frac{\partial \widehat{f}(q, \sigma)}{\partial q}dq}
{|p|^{2s_{1}}+|p|^{2s_{2}}}\chi_{\{|p|\leq 1\}}\Bigg|\leq
\frac{\|xf(x)\|_{L^{1}({\mathbb R}^{3})}}{(2\pi)^{\frac{3}{2}}}|p|^{1-2s_{1}}
\chi_{\{|p|\leq 1\}}\in L^{2}({\mathbb R}^{3}).
$$
It remains to analyze the term
\begin{equation}
\label{fh0s12}
\frac{\widehat{f}(0)}{|p|^{2s_{1}}+|p|^{2s_{2}}}\chi_{\{|p|\leq 1\}}.
\end{equation}
It can be easily verified that (\ref{fh0s12}) is contained in
$L^{2}({\mathbb R}^{3})$ if and only if $\widehat{f}(0)=0$. This is equivalent to
orthogonality condition (\ref{oc1}). \hfill\lanbox

\bigskip

Note that by proving the lemma above we establish the solvability of
equation (\ref{lp}) in $H^{2s_{2}}({\mathbb R}^{3})$
for all the values of the powers of the fractional Laplacians
$\displaystyle{0<s_{1}<s_{2}<1}$, so that no orthogonality relations are required
for the right side $f(x)$ when
$\displaystyle{s_{1}\in \Big(0, \frac{3}{4}\Big)}$. But if
$\displaystyle{s_{1}\in \Big[\frac{3}{4}, 1\Big)}$, the argument of the proof
relies on the single orthogonality condition (\ref{oc1}).
This is similar to the case when the Poisson
type equation is considered with a single fractional Laplacian in
${\mathbb R}^{3}$ (see Theorem 1.1 of ~\cite{VV19}, also ~\cite{VV15}).
For the solvability of the problem similar to (\ref{lp}) containing a scalar
potential see ~\cite{EV22}, also ~\cite{V2024}.

\medskip

Let us write down the corresponding sequence of the approximate equations
related to problem (\ref{lp}) with $n\in {\mathbb N}$ as
\begin{equation}
\label{lpn}
[{(-\Delta)}^{s_{1}}+{(-\Delta)}^{s_{2}}]u_{n}=f_{n}(x), \quad
x\in {\mathbb R}^{3}, \quad 0<s_{1}<s_{2}<1.
\end{equation}
The right sides of (\ref{lpn}) tend to the right side of (\ref{lp}) in
$L^{2}({\mathbb R}^{3})$ as
$n\to \infty$. We demonstrate that under the certain technical assumptions
each equation (\ref{lpn}) has a unique solution
$u_{n}(x)\in H^{2s_{2}}({\mathbb R}^{3})$, limiting problem (\ref{lp}) admits
a unique solution $u(x)\in H^{2s_{2}}({\mathbb R}^{3})$ and
$u_{n}(x)\to u(x)$ in $H^{2s_{2}}({\mathbb R}^{3})$ as $n\to \infty$. This is
the so called {\it solvability in the sense of sequences} for equation
(\ref{lp}). The final statement of our work is as follows.

\bigskip

\noindent
{\bf Lemma 4.2.} {\it  Let  $n\in {\mathbb N}, \ 0<s_{1}<s_{2}<1, \ f_{n}(x):
{\mathbb R}^{3}\to {\mathbb R}, \ f_{n}(x)\in L^{2} ({\mathbb R}^{3})$,
such that $f_{n}(x)\to f(x)$ in $L^{2} ({\mathbb R}^{3})$ as $n\to \infty$.

\medskip

\noindent
a) Suppose that $\displaystyle{s_{1}\in \Big(0, \frac{3}{4}\Big)}$ and in
addition
$f_{n}(x)\in L^{1} ({\mathbb R}^{3}), \ n\in {\mathbb N}$, so that
$f_{n}(x)\to f(x)$ in $L^{1} ({\mathbb R}^{3})$ as $n\to \infty$. Then
problems (\ref{lp}) and (\ref{lpn}) possess unique solutions
$u(x)\in H^{2s_{2}} ({\mathbb R}^{3})$ and $u_{n}(x)\in H^{2s_{2}} ({\mathbb R}^{3})$
respectively, such that $u_{n}(x)\to u(x)$ in $H^{2s_{2}} ({\mathbb R}^{3})$
as $n\to \infty$.

\medskip

\noindent
b) Suppose that
$\displaystyle{s_{1}\in \Big[\frac{3}{4}, 1 \Big)}$ and additionally
$xf_{n}(x)\in L^{1} ({\mathbb R}^{3}), \ n\in {\mathbb N}$, so that
$xf_{n}(x)\to xf(x)$ in $L^{1} ({\mathbb R}^{3})$
as $n\to \infty$ and the orthogonality relations
\begin{equation}
\label{oc2}
(f_{n}(x), 1)_{L^{2} ({\mathbb R}^{3})}=0
\end{equation}
are valid for all $n\in {\mathbb N}$. Then
equations (\ref{lp}) and (\ref{lpn}) admit unique solutions
$u(x)\in H^{2s_{2}} ({\mathbb R}^{3})$ and $u_{n}(x)\in H^{2s_{2}} ({\mathbb R}^{3})$
respectively, so that $u_{n}(x)\to u(x)$ in $H^{2s_{2}} ({\mathbb R}^{3})$
as $n\to \infty$.}

\bigskip

\noindent
{\it Proof.} Let us verify that if problems
(\ref{lp}) and (\ref{lpn}) have unique solutions
$u(x)\in H^{2s_{2}}({\mathbb R}^{3})$ and
$u_{n}(x)\in H^{2s_{2}}({\mathbb R}^{3}), \ n\in {\mathbb N}$ respectively and
$u_{n}(x)\to u(x)$ in $L^{2}({\mathbb R}^{3})$ as
$n\to \infty$, then we have $u_{n}(x)\to u(x)$ in
$H^{2s_{2}}({\mathbb R}^{3})$ as $n\to \infty$ as well. Clearly, by means of
(\ref{lpn}) and (\ref{lp})
$$
[(-\Delta)^{s_{1}}+(-\Delta)^{s_{2}}](u_{n}(x)-u(x))=f_{n}(x)-f(x).
$$
By applying the standard Fourier transform (\ref{f}), we obtain that
$$
\|(-\Delta)^{s_{2}}(u_{n}(x)-u(x))\|_{L^{2}({\mathbb R}^{3})}\leq
\|f_{n}(x)-f(x)\|_{L^{2}({\mathbb R}^{3})}\to 0, \quad n\to \infty
$$
as assumed. Norm definition (\ref{n}) implies
$u_{n}(x)\to u(x)$ in $H^{2s_{2}}({\mathbb R}^{3})$ as $n\to \infty$.

Let us apply (\ref{f}) to both sides of equations (\ref{lp}) and (\ref{lpn})
to derive
\begin{equation}
\label{unuhp}
\widehat{u_{n}}(p)-\widehat{u}(p)=\frac{\widehat{f_{n}}(p)-\widehat{f}(p)}
{|p|^{2s_{1}}+|p|^{2s_{2}}}\chi_{\{|p|\leq 1\}}+\frac{\widehat{f_{n}}(p)-\widehat{f}(p)}
{|p|^{2s_{1}}+|p|^{2s_{2}}}\chi_{\{|p|>1\}}.
\end{equation}
Obviously, the second term in the right side of (\ref{unuhp}) can be bounded
from above in the absolute value by
$\displaystyle{\frac{|\widehat{f_{n}}(p)-\widehat{f}(p)|}{2}}$, such that
$$
\Bigg\|\frac{\widehat{f_{n}}(p)-\widehat{f}(p)}{|p|^{2s_{1}}+|p|^{2s_{2}}}
\chi_{\{|p|>1\}}\Bigg\|_{L^{2}({\mathbb R}^{3})}\leq \frac{1}{2}
\|f_{n}(x)-f(x)\|_{L^{2}({\mathbb R}^{3})}\to 0, \quad n\to \infty
$$
as we assume.

First we consider the case when
$\displaystyle{s_{1}\in \Big(0, \frac{3}{4}\Big)}$. According to the result
of Lemma 4.1 above, equations (\ref{lp}) and (\ref{lpn}) admit unique solutions
$u(x)\in H^{2s_{2}}({\mathbb R}^{3})$ and
$u_{n}(x)\in H^{2s_{2}}({\mathbb R}^{3}), \ n\in {\mathbb N}$ respectively.

Evidently, the first term in the right side of (\ref{unuhp}) can be estimated
from above in the absolute value via (\ref{fub}) by
$$
\frac{\|f_{n}(x)-f(x)\|_{L^{1}({\mathbb R}^{3})}}{(2\pi)^{\frac{3}{2}}|p|^{2s_{1}}}
\chi_{\{|p|\leq 1\}},
$$
so that
$$
\Bigg\|\frac{\widehat{f_{n}}(p)-\widehat{f}(p)}{|p|^{2s_{1}}+|p|^{2s_{2}}}
\chi_{\{|p|\leq 1\}}\Bigg\|_{L^{2}({\mathbb R}^{3})}\leq
\frac{\|f_{n}(x)-f(x)\|_{L^{1}({\mathbb R}^{3})}}{\sqrt{2}\pi \sqrt{3-4s_{1}}}
\to 0, \quad n\to \infty
$$
due to the one of our assumptions. This means that
$$
u_{n}(x)\to u(x) \quad in \quad L^{2}({\mathbb R}^{3}) \quad as \quad
n\to \infty
$$
in the situation when
$\displaystyle{s_{1}\in \Big(0, \frac{3}{4}\Big)}$. By virtue of the reasoning
above, this completes the proof of the case a) of the lemma.

Let us turn our attention to the situation when
$\displaystyle{s_{1}\in \Big[\frac{3}{4}, 1\Big)}$. Recall parts a) and b) of
Lemma 4.1 of ~\cite{VV19}. Under the given conditions, we have
$f_{n}(x)\in L^{1}({\mathbb R}^{3}), \ n\in {\mathbb N}$ and
$f_{n}(x)\to f(x)$ in $L^{1}({\mathbb R}^{3})$ as $n\to \infty$. Using
(\ref{oc2}), we obtain
$$
|(f(x), 1)_{L^{2}({\mathbb R}^{3})}|=|(f(x)-f_{n}(x), 1)_{L^{2}({\mathbb R}^{3})}|\leq
\|f_{n}(x)-f(x)\|_{L^{1}({\mathbb R}^{3})}\to 0
$$
as $n\to \infty$. Thus,
\begin{equation}
\label{ocl}
(f(x), 1)_{L^{2}({\mathbb R}^{3})}=0
\end{equation}
holds. By means of the result of part b) of Lemma 4.1, problems (\ref{lp}) and
(\ref{lpn}) possess unique solutions
$u(x)\in H^{2s_{2}}({\mathbb R}^{3})$ and
$u_{n}(x)\in H^{2s_{2}}({\mathbb R}^{3}), \ n\in {\mathbb N}$ respectively.

Orthogonality relations (\ref{ocl}) and (\ref{oc2}) imply that
$$
\widehat{f}(0)=0, \quad \widehat{f_{n}}(0)=0, \quad n\in {\mathbb N}.
$$
This enables us to express
\begin{equation}
\label{fhfnh}
\widehat{f}(p)=\int_{0}^{|p|}\frac{{\partial \widehat{f}}(q, \sigma)}
{\partial q}dq, \quad
\widehat{f_{n}}(p)=\int_{0}^{|p|}\frac{{\partial \widehat{f_{n}}}(q, \sigma)}
{\partial q}dq, \quad n\in {\mathbb N},
\end{equation}
so that the first term in the right side of (\ref{unuhp}) is given by
\begin{equation}
\label{intdfndf}
\frac
{\int_{0}^{|p|}\Big[\frac{{\partial \widehat{f_{n}}}(q, \sigma)}{\partial q}-
\frac{{\partial \widehat{f}}(q, \sigma)}{\partial q}\Big]dq}
{|p|^{2s_{1}}+|p|^{2s_{2}}}\chi_{\{|p|\leq 1\}}.
\end{equation}
By virtue of the definition of the standard Fourier transform (\ref{f}), we
easily derive that
\begin{equation}
\label{dfnfhp}
\Bigg|\frac{\partial \widehat{f_{n}}(p)}{\partial |p|}-
\frac{\partial \widehat{f}(p)}{\partial |p|}\Bigg|\leq
\frac{\|xf_{n}(x)-xf(x)\|_{L^{1}({\mathbb R}^{3})}}
{(2\pi)^{\frac{3}{2}}}.
\end{equation}
This implies
$$
\Bigg|\frac
{\int_{0}^{|p|}\Big[\frac{{\partial \widehat{f_{n}}}(q, \sigma)}{\partial q}-
\frac{{\partial \widehat{f}}(q, \sigma)}{\partial q}\Big]dq}
{|p|^{2s_{1}}+|p|^{2s_{2}}}\chi_{\{|p|\leq 1\}}\Bigg|\leq
\frac{\|xf_{n}(x)-xf(x)\|_{L^{1}({\mathbb R}^{3})}}
{(2\pi)^{\frac{3}{2}}}|p|^{1-2s_{1}}\chi_{\{|p|\leq 1\}},
$$
so that
$$
\Bigg\|\frac
{\int_{0}^{|p|}\Big[\frac{{\partial \widehat{f_{n}}}(q, \sigma)}{\partial q}-
\frac{{\partial \widehat{f}}(q, \sigma)}{\partial q}\Big]dq}
{|p|^{2s_{1}}+|p|^{2s_{2}}}\chi_{\{|p|\leq 1\}}\Bigg\|_{L^{2}({\mathbb R}^{3})}\leq
\frac{\|xf_{n}(x)-xf(x)\|_{L^{1}({\mathbb R}^{3})}}{\pi \sqrt{2(5-4s_{1})}}\to 0
$$
as $n\to \infty$ via the one of our assumptions. Therefore,
$$
u_{n}(x)\to u(x) \quad in \quad L^{2}({\mathbb R}^{3}) \quad as \quad n\to \infty
$$
in the case when
$\displaystyle{s_{1}\in \Big[\frac{3}{4}, 1\Big)}$. By means of the argument
above, this completes the proof of the part b) of our lemma.
\hfill\lanbox

\bigskip


\section*{Acknowledgements}

The first author is grateful to Israel Michael Sigal for the partial support
by the NSERC grant NA 7901.
The second author has been supported by the RUDN University Strategic Academic
Leadership Program.

\end{document}